\documentstyle[12pt, amsfonts]{amsart}
\begin{document}
\def\R{{\mathbb R}}
\def\Z{{\mathbb Z}}
\def\C{{\mathbb C}}
\newcommand{\trace}{\rm trace}
\newcommand{\Ex}{{\mathbb{E}}}
\newcommand{\Prob}{{\mathbb{P}}}
\newcommand{\E}{{\cal E}}
\newcommand{\F}{{\cal F}}
\newtheorem{df}{Definition}
\newtheorem{theorem}{Theorem}
\newtheorem{lemma}{Lemma}
\newtheorem{pr}{Proposition}
\newtheorem{co}{Corollary}
\def\n{\nu}
\def\sign{\mbox{ sign }}
\def\a{\alpha}
\def\N{{\mathbb N}}
\def\A{{\cal A}}
\def\L{{\cal L}}
\def\X{{\cal X}}
\def\F{{\cal F}}
\def\c{\bar{c}}
\def\v{\nu}
\def\d{\delta}
\def\diam{\mbox{\rm dim}}
\def\vol{\mbox{\rm Vol}}
\def\b{\beta}
\def\t{\theta}
\def\l{\lambda}
\def\e{\varepsilon}
\def\colon{{:}\;}
\def\pf{\noindent {\bf Proof :  \  }}
\def\endpf{ \begin{flushright}
$ \Box $ \\
\end{flushright}}

\title[Slicing inequalities for subspaces of $L_p$]
{Slicing inequalities for subspaces of $L_p$}

\author{Alexander Koldobsky}

\address{Department of Mathematics\\ 
University of Missouri\\
Columbia, MO 65211}

\email{koldobskiya@@missouri.edu}

\begin{abstract}  We prove slicing inequalities for measures of the unit balls of subspaces of $L_p$,
$-\infty<p<\infty.$
For example, for every $k\in \N$ there exists a constant $C(k)$ such that for every $n\in \N,\ k<n$,
every convex $k$-intersection body (unit ball of a normed subspace of $L_{-k})$
$L$ in $\R^n$ and every measure $\mu$ with non-negative even continuous
density in $\R^n,$
$$\mu(L)\ \le\ C(k) \max_{\xi \in S^{n-1}} 
\mu(L\cap \xi^\bot)\ |L|^{1/n} \ ,$$
where  $\xi^\bot$ is the central hyperplane in $\R^n$ perpendicular to $\xi,$ and
$|L|$ is the volume of $L.$ This and other results are in the spirit of the hyperplane problem of Bourgain.
The proofs are based on stability inequalities for intersection bodies and estimates for the Banach-Mazur
distance from the unit ball of a subspace of $L_p$ to the class of intersection bodies. 
\end{abstract}  
\maketitle

\section{Introduction}
The hyperplane problem \cite{Bo1, Bo2, Ba1, MP} asks whether there exists 
an absolute constant $C$ so that for any origin-symmetric convex body $K$ in $\R^n$
\begin{equation} \label{hyper}
|K|^{\frac {n-1}n} \le C \max_{\xi \in S^{n-1}} |K\cap \xi^\bot|,
\end{equation}
where  $\xi^\bot$ is the central hyperplane in $\R^n$ perpendicular to $\xi,$ and
$|K|$ stands for volume of proper dimension.
The problem is still open. Bourgain \cite{Bo3} proved that $C\le O(n^{1/4}\log n),$ and Klartag \cite{Kl}
removed the  logarithmic term and established the best-to-date estimate $C\le  O(n^{1/4}).$ 
We refer the reader to [BGVV] for the history and partial results.

Inequality (\ref{hyper}) has been proved for several special classes of bodies. 
For instance, if $K$ is an intersection body,
the inequality is known with the best possible constant (see \cite[Th. 9.4.11]{G}):
\begin{equation} \label{hyper-int}
|K|^{\frac {n-1}n} \le c_n \max_{\xi \in S^{n-1}} |K\cap \xi^\bot|,
\end{equation}
where $c_n= |B_2^n|^{\frac {n-1}n} / |B_2^{n-1}|<1,$ and $B_2^n$ is the unit Euclidean ball.
One possible path to attack the hyperplane problem is through the scale of finite dimensional
subspaces of $L_p,$ because one approximately gets all origin-symmetric convex bodies as $p\to +\infty.$ 
It was shown by Ball \cite{Ba2} that the hyperplane problem has an affirmative answer for the unit
balls of finite dimensional subspaces 
of $L_p,\ 1\le p \le 2.$ Junge \cite{J} extended this result to every $p\in (1,\infty)$ with the constant $C$
depending on $p$ and going to infinity when $p\to \infty.$ Milman \cite{Mi1} gave a different proof for
subspaces of $L_p,\ 2\le p <\infty$ with the constant $C\le O(\sqrt{p}).$ Another proof of this estimate
can be found in \cite{KPY}.

The paper \cite{KPY} also explored a different path, through $n$-dimensional
spaces that embed in $L_{-p},\ 0<p<n.$ The latter concept was introduced in \cite{K2,K4}; see definition
and basic facts below. In particular, a normed space embeds in $L_{-k}$ iff its unit ball belongs
to the class of $k$-intersection bodies introduced in \cite{K3,K4} in connection with generalizations of the Busemann-Petty
problem. The path through $L_{-p}$ presents some hope because every $n$-dimensional
normed space embeds in $L_{-p}$ with $n-3 \le  p < n.$ However, the paper \cite{KPY} did not go very
far, as the hyperplane conjecture was confirmed there only for the unit balls of spaces that embed
in $L_{-p},\ 0<p<2,$ with the constant going to infinity as $p\to -2.$

In this note we prove the hyperplane conjecture for the unit balls of normed spaces that embed in
$L_{-p}$ for arbitrary $p>0,$ and, in particular, for $k$-intersection bodies for all $k\in \N.$ Moreover, we prove the result
with arbitrary measures in place of volume. 

The first result on the hyperplane problem for arbitrary measures was proved in \cite{K6}.
If  $K$ is an intersection body in $\R^n$  and $\mu$ is any measure with even continuous 
density $f$ in $\R^n$ (so that $\mu(B)=\int_B f$ for every compact set $B\subset \R^n),$ then
\begin{equation} \label{arbmeas}
\mu(K) \le \frac n{n-1} c_n \max_{\xi \in S^{n-1}} \mu(K\cap \xi^\bot)\ |K|^{1/n},
\end{equation}
where
$c_n$ is the same constant as in (\ref{hyper-int}). The constant in (\ref{arbmeas}) is sharp
and differs from the case of volume in (\ref{hyper-int}) by a small factor of $n/(n-1).$ Recall
that $c_n<1.$

The condition that $K$ is an intersection body in inequality (\ref{arbmeas}) was removed in \cite{K7}
at the expense of extra $\sqrt{n}$ in the right-hand side. For any origin-symmetric convex body $K$
in $\R^n$ and any measure $\mu$ with even continuous density in $\R^n$
\begin{equation} \label{sqrtn2}
\mu(K)\ \le\ \sqrt{n} \frac n{n-1} c_n\max_{\xi \in S^{n-1}} 
\mu(K\cap \xi^\bot)\ |K|^{1/n}.
\end{equation}
Analogs of (\ref{sqrtn2}) for sections of lower dimensions and for complex convex bodies 
were proved in \cite{K8}.

Our main result is the following inequality. For any $p>0$ there exists a constant $C(p)$ such that for any $n\in \N,\ n>p,$
any convex body $L$ in $\R^n$ that is the unit ball of a normed space embedding in $L_{-p},$ any $1\le m \le n-1,$ and
any measure $\mu$ with even continuous non-negative density in $\R^n,$
\begin{equation} \label{slicing-p}
\mu(L)\ \le\ (C(p))^m \frac n{n-m} \max_{H\in Gr_{n-m}} 
\mu(L\cap H)\ |L|^{m/n},
\end{equation}
where $Gr_{n-m}$ is the Grassmanian of $(n-m)$-dimensional subspaces of $\R^n.$
In the case where $p=k\in \N$ inequality (\ref{slicing-p}) applies to any convex $k$-intersection body $L.$

We also show that the constant $\sqrt{n}$ in (\ref{sqrtn2}) can be replaced by $n^{1/2-1/p}$ when
$K$ is the unit ball of an $n$-dimensional subspace of $L_p, p>2.$ 
The proofs are based on a stability result for intersection bodies, which was proved in \cite{K8}
and continued the study of stability in volume comparison problems initiated in \cite{K5}; see also \cite{K9}.

\section{Embedding in $L_{-p}$ and $k$-intersection bodies}

We need several definitions and facts.
A closed bounded set $K$ in $\R^n$ is called a {\it star body}  if 
every straight line passing through the origin crosses the boundary of $K$ 
at exactly two points different from the origin, the origin is an interior point of $K,$
and the {\it Minkowski functional} 
of $K$ defined by 
$$\|x\|_K = \min\{a\ge 0:\ x\in aK\}$$
is a continuous function on $\R^n.$ 

The {\it radial function} of a star body $K$ is defined by
$$\rho_K(x) = \|x\|_K^{-1}, \qquad x\in \R^n,\ x\neq 0.$$
If $x\in S^{n-1}$ then $\rho_K(x)$ is the radius of $K$ in the
direction of $x.$

The class of intersection bodies was introduced by Lutwak \cite{L}.
Let $K, L$ be origin-symmetric star bodies in $\R^n.$ We say that $K$ is the 
intersection body of $L$ if the radius of $K$ in every direction is 
equal to the $(n-1)$-dimensional volume of the section of $L$ by the central
hyperplane orthogonal to this direction, i.e. for every $\xi\in S^{n-1},$
\begin{equation} \label{intbodyofstar}
\rho_K(\xi)= \|\xi\|_K^{-1} = |L\cap \xi^\bot|.
\end{equation} 
All bodies $K$ that appear as intersection bodies of different star bodies
form {\it the class of intersection bodies of star bodies}. The class of {\it intersection bodies} 
can be defined as the closure of the class of intersection bodies of star bodies
in the radial metric 
$$\rho(K,L)=\sup_{\xi\in S^{n-1}} \left|\rho_K(\xi)-\rho_L(\xi)\right|.$$  

A more general concept of a $k$-intersection body was introduced 
in \cite{K3,K4}. For an integer $k,\ 1\le k <n$ and star bodies $D,L$ in $\R^n,$
we say that $D$ is the $k$-intersection body of $L$ if for every $(n-k)$-dimensional
subspace $H$ of $\R^n,$ 
$$|D\cap H^\bot|= |L\cap H|,$$
where $H^\bot$ is the $k$-dimensional subspace orthogonal to $H.$
Taking the closure in the radial metric of the class of all $D$'s that appear
as $k$-intersection bodies of star bodies, we define the class of 
{\it $k$-intersection bodies} (the original definition in \cite{K3,K4} was
different; the equivalence of definitions was proved by Milman \cite{Mi2}). 
If $k=1$ one gets the usual intersection bodies. 
Intersection bodies played a crucial role in the solution of the
Busemann-Petty problem and its generalizations; see \cite[Ch. 5]{K1}.

Another generalization of intersection bodies was introduced by Zhang \cite{Z} in connection
with the lower dimensional Busemann-Petty problem.
For $1\le m \le n-1,$ denote by $Gr_{n-m}$ the Grassmanian of $(n-m)$-dimensional
subspaces of $\R^n.$ The {\it $(n-m)$-dimensional spherical Radon transform} 
$R_{n-m}:C(S^{n-1})\mapsto C(Gr_{n-m})$  
is a linear operator defined by
$$R_{n-m}g (H)=\int_{S^{n-1}\cap H} g(x)\ dx,\quad \forall  H\in Gr_{n-m}$$
for every function $g\in C(S^{n-1}).$ Denote by
$$R_{n-m}\left(C(S^{n-1})\right)=X\subset C(Gr_{n-m}).$$
Let $M^+(X)$ be the space of linear positive continuous functionals on $X$, i.e. for every
$\nu\in M^+(X)$ and non-negative function $f\in X$, we have $\nu(f)\geq0$.

An origin-symmetric star body $K$ in $\R^n$ is called a {\it generalized $m$-intersection
body} if there exists a functional $\nu\in M^+(X)$ so that for every $g\in C(S^{n-1})$,
\begin{equation} \label{defintbody}
\int_{S^{n-1}} \|x\|_K^{-m} g(x)\ dx=\nu(R_{n-m}g).
\end{equation}
It was shown in \cite{K4} (see \cite[p. 92]{K1} or \cite{Mi2}
for different proofs) that every generalized $m$-intersection body is 
an $m$-intersection body. Note that generalized 1-intersection bodies 
coincide with the original intersection bodies of Lutwak. We denote the 
class of generalized $m$-intersection bodies by $BP_{m}.$

The concept of embedding of finite dimensional normed spaces in $L_{-p}$ with $p>0$
was introduced in \cite{K2,K4} as an analytic 
extension of embedding of normed spaces into $L_p$ with $p>0.$ It is a well-known fact
going back to P.Levy (see for example \cite[Lemma 6.4]{K1}) that an $n$-dimensional normed space 
$(\R^n,\|\cdot\|)$ embeds 
in $L_p([0,1]),\ p>-1$ if and only if there exists a finite Borel measure $\mu$ on the 
sphere $S^{n-1}$ in $\R^n$ so that for every $x\in \R^n,$
\begin{equation} \label{levy}
\|x\|^p =\int_{S^{n-1}} |(x,\xi)|^p d\mu(\xi).
\end{equation}
We say that a normed space embeds in $L_p$ if its norm satisfies (\ref{levy}),
without associating $L_p$ with a measure space.
Embedding in $L_{-p}$ is defined as an analytic extension of (\ref{levy}).

\begin{df} Let $X$ be an $n$-dimensional normed space, and $0<p<n.$ 
We say that $X$ embeds in $L_{-p}$
if there exists a finite Borel measure $\mu $ on $S^{n-1}$ so that, for every even Schwartz test 
function $\phi\in {\mathcal S}(\R^n),$
$$
\int_{\R^n} \|x\|_X^{-p} \phi(x)\ dx = \int_{S^{n-1}}  \left(
\int_{\R} |t|^{p-1} \hat{\phi}(t\theta)\ dt \right) d\mu(\theta),
$$
where $\hat{\phi}$ stands for the Fourier transform.
\end{df}
   
A connection between
intersection bodies and embedding in $L_p$ was established in \cite[Th. 3]{K4}.

\begin{pr} \label{connection} Let $1\le k < n.$ The following are equivalent:
\item{(i)} An origin symmetric star body $D$ in $\R^n$ is a $k$-intersection
body;
\item{(ii)} $\|\cdot\|_D^{-k}$ represents a positive definite distribution;
\item{(iii)} The space $(\R^n,\|\cdot\|_D)$ embeds in $L_{-k}.$
\end{pr} 

The advantage of this connection (and, consequently, of introducing embeddings
in $L_{-p})$ is that now one can try to extend to negative values of $p$ 
different results about usual $L_p$-spaces. Every such extension gives new information about 
intersection bodies. Let us show several examples of this approach.

A well-known simple fact is that every two-dimensional normed space embeds in $L_1.$
It can be proved (see \cite[Theorem 4.13]{K1}) 
that for every symmetric convex body $K$ in $\R^n$
and every $p\in [n-3,n),$ the function $\|\cdot\|_K^{-p}$ represents a
positive definite distribution, so by Theorem \ref{connection} every $n$-dimensional
Banach space embeds in $L_{-n+3}.$ Putting $n=2$ we get 
the property of two-dimensional spaces mentioned above.
Putting $n=4$ we see that every four-dimensional normed space embeds
in $L_{-1}.$ By Theorem \ref{connection}, every four-dimensional symmetric convex body
is a 1-intersection body, which solves in affirmative the 
critical four-dimensional case of the Busemann-Petty problem; see \cite[Ch. 5]{K1}. 

Another well-known property of $L_p$-spaces is that, for any $0<p<q\le 2,$ 
the space $L_q$ embeds isometrically in $L_p,$ so $L_p$-spaces become larger
when $p$ decreases from 2. 
This result can be extended to negative $p$ as follows.
\begin{pr} \label{LpLq} (\cite[Th. 6.17]{K1}) Every $n$-dimensional subspace of $L_q,\ 0<q\le 2$
embeds in $L_{-p}$ for every $p\in (0,n).$ Hence, the unit ball of
every $n$-dimensional subspace of $L_q,\ 0<q\le 2$ is a $k$-intersection
body for every $1\le k < n.$ 
\end{pr}

The factorization theorem 
of Maurey and Nikishin \cite{Ma, Ni} implies that, 
for $0<p<q<1,$ every Banach subspace of $L_p$ is isomorphic 
to a subspace of $L_q$. The following extension
to negative $p$ was proved in \cite{KK} (see also \cite[Section 6.3]{K1}).  
Recall that the geometric distance between two star bodies $K$ and $L$, $d_G(K,L),$ is defined as the infimum
of positive numbers $r$ such that there exists some $a >0$ so that
$$K \subset aL \subset rK,$$ and
the Banach-Mazur distance is defined by
$$d_{BM}(K,L) = \inf_{T \in GL_n}
d_G(K,TL).$$

\begin{pr}\label{factor2} (\cite{KK}) If $-\infty<p<q<1$, $p\not=0$ and $q>0$ then there exists a     
constant $C=C(p,q)$ so that, for any $n$ with $-n<p,$ if $X$ is an $n$-dimensional normed space that 
embeds into $L_p,$ then there is a subspace $Y$ of $L_q$ so that $d_{BM}(K,L)\le C,$ where
$K$ and $L$ are the unit balls of $X$ and $Y,$ respectively.      
\end{pr}   

We refer the reader to \cite{K1, KY, Mi2, Mi3, Y1, Y2, KZ, Z}
for more about $k$-intersection bodies, embedding in $L_{-p}$ and their applications.

\section{Slicing inequalities}
 
We use  the following stability result for generalized intersection bodies proved in \cite{K8}.  
\begin{pr}\label{stab1}{\bf (\cite{K8})}
Suppose that $1\le m \le n-1,$ $K$ is a generalized $m$-intersection body in $\R^n,$  $f$
is an even continuous function on $K,$ $f\ge 1$ everywhere on $K,$ and $\e>0.$ If
$$
\int_{K\cap H} f \ \le\ |K\cap H| +\e,\qquad \forall H\in Gr_{n-m},
$$
then
$$
\int_K f\ \le\ |K| + \frac {n}{n-m}\ c_{n,m}\ |K|^{m/n}\e,
$$
where $c_{n,m}= |B_2^n|^{\frac {n-m}n}/|B_2^{n-m}|;$ note that $c_{n,m}\in ((\frac 1{e})^{m/2}, 1)$
and that $c_{n,1}$ equals to the constant $c_n$ from (\ref{hyper-int}). 
\end{pr}

We define the Banach-Mazur distance from a star body $L$ in $\R^n$ to the class $BP_{m}$ of
generalized $m$-intersection bodies by
$$d_{BM}(L,BP_{m})= \inf_{K\in BP_{m}} d_{BM}(L,K).$$

\begin{theorem}\label{bm} Suppose that $1\le m \le n-1$ and  $L$ is an origin-symmetric star body in $\R^n$
such that $d_{BM}(L, BP_{m})< d$ for some $d>0.$ Then for any measure $\mu$ with even continuous non-negative density on $L$
\begin{equation} \label{sqrtn}
\mu(L)\ \le\  d^m\ \frac n{n-m} c_{n,m} \max_{H\in Gr_{n-m}}  \mu(L\cap H)\ |L|^{m/n}.
\end{equation}
\end{theorem}

\pf  Since any linear transformation of a generalized $m$-intersection body is again a generalized $m$-intersection body
(see \cite[Corollary 6.8]{GZ}),
there exists $K\in BP_{m}$ such that
$$\frac 1d \ K \subset L \subset K.$$ 
Let $g$ be the density of $\mu,$ and let $f= \chi_K + g \chi_L,$ where $\chi_K,\ \chi_L$ are the indicator functions of $K$ and $L.$ 
Then $f\ge 1$ everywhere on $K.$ Put 
$$\e=\max_{H\in Gr_{n-m}} \left(\int_{K\cap H} f - |K\cap H| \right)= \max_{H\in Gr_{n-m}} \int_{L\cap H} g.$$
Now we can apply Proposition \ref{stab1} to $f,K,\e$ (the function $f$ is not necessarily continuous on $K,$ 
but the result holds by a simple approximation argument). We get
$$\mu(L)= \int_L g = \int_K f  -\ |K|$$
$$ \le \frac n{n-m} c_{n,m}  |K|^{m/n}\max_{H\in Gr_{n-m}} \int_{L\cap H} g$$
$$ \le d^m\ \frac n{n-m} c_{n,m} |L|^{m/n}\max_{H\in Gr_{n-m}} \mu(L\cap H),$$
because $K\subset d L,$ so $|K|\le d^n  |L|.$ \qed 

Combining Theorem \ref{bm} with Proposition \ref{factor2}, we prove a slicing inequality for
subspaces of $L_{-p}.$ 

\begin{theorem} \label{main-lp} For any $p>0$ there exists a constant $C(p)$ such that for any $n\in \N, n>p,$
any convex body $L$ in $\R^n$ that is the unit ball of a normed space embedding in $L_{-p},$ any $1\le m \le n-1,$ and
any measure $\mu$ with even continuous non-negative density in $\R^n,$
\begin{equation} \label{slicing-negp}
\mu(L)\ \le\ (C(p))^m \frac n{n-m} c_{n,m} \max_{H\in Gr_{n-m}} 
\mu(L\cap H)\ |L|^{m/n}.
\end{equation}
Recall that $c_{n,m}<1.$
\end{theorem}

\pf By Proposition \ref{factor2}, there exists a body $K$ which is the unit ball of a normed space 
that embeds in $L_{1/2}$ and such that $d_{BM}(K,L)< C(p,1/2) =: C(p).$ By Proposition \ref{LpLq}, $K$ 
is an intersection body. By \cite[Lemma 6.1]{GZ}, every intersection body is a generalized 
$m$-intersection body for every $1\le m \le n-1,$ so $K\in BP_{m}.$ The result follows from Theorem \ref{bm}. \qed

The author is unable to estimate the behavior of the constant $C(p)$ as $p\to \infty,$ because the proof 
of Proposition \ref{factor2} is rather involved.
For $p=k\in \N,$ Theorem \ref{main-lp} can be reformulated using Proposition \ref{connection}.
\begin{co} \label{main-kint} For any $k\in \N$ there exists a constant $C(k)$ such that for any $n\in \N,\ n>k,$ 
any convex $k$-intersection body $L$ in $\R^n,$ any $1\le m \le n-1,$ and any measure $\mu$ with even continuous density
in $\R^n$
\begin{equation} \label{slicing-kint}
\mu(L)\ \le\ (C(k))^m \frac n{n-m} c_{n,m} \max_{H\in Gr_{n-m}} 
\mu(L\cap H)\ |L|^{m/n}.
\end{equation}
\end{co}

Let us formulate the hyperplane cases ($m=1$) of Theorem \ref{main-lp} and Corollary \ref{main-kint} separately,
as they relate to the hyperplane problem.
\begin{co} \label{main-lp} For any $p>0$ there exists a constant $D(p)$ such that for any $n\in \N,\ n>p,$
any convex body $L$ in $\R^n$ that is the unit ball of a normed space embedding in $L_{-p}$ and
any measure $\mu$ with even continuous non-negative density in $\R^n,$
\begin{equation} \label{slicing-negp}
\mu(L)\ \le\ D(p) \max_{\xi \in S^{n-1}} 
\mu(L\cap \xi^\bot)\ |L|^{1/n}.
\end{equation}
When $p=k\in \N,\ 1\le k \le n-1$ the latter inequality applies to any convex $k$-intersection body $L$ in $\R^n;$
see the formulation in the abstract.
\end{co}

Finally, let us prove a slicing inequality for subspaces of $L_p,\ p>2$ with arbitrary measures.
\begin{theorem} \label{p>2} Suppose that $p>2,$ $L$ is the unit ball of an $n$-dimensional subspace of $L_p,$
$1\le m \le n-1,$ and $\mu$ is a measure with even continuous density in $\R^n.$ Then
$$\mu(L) \le n^{m/2-m/p} \frac n{n-m} c_{n,m} \max_{H\in Gr_{n-m}} \mu(L\cap H) |L|^{m/n}.$$
\end{theorem}

\pf By a result of Lewis \cite{Le} (see also \cite{SZ} for a different proof), $d_{BM}(L, B_2^n) \le n^{1/2-1/p}.$
Since $B_2^n\in BP_m$ for every $m$, the result follows from Theorem \ref{bm}.
\qed
\bigbreak
The hyperplane case is as follows.
\begin{co} If $L$ is the unit ball of an $n$-dimensional subspace of $L_p,\ p>2,$
then for any measure $\mu$ with even continuous density in $\R^n$
$$\mu(L)\le n^{1/2-1/p} \frac n{n-1} c_n \max_{\xi\in S^{n-1}} \mu(L\cap \xi^\bot)\ |L|^{1/n}.$$ 
\end{co}
\bigbreak

\noindent{\bf Remark.} In Theorems \ref{main-lp} and \ref{p>2} we consider $L_p$-spaces with $p<0$
and $p>2.$ As follows from Proposition \ref{LpLq}, the unit balls of subspaces of $L_p$ with $0<p\le 2$ are 
intersection bodies and satisfy inequalities (\ref{hyper-int}) and (\ref{arbmeas}). In the remaining case $p=0,$ the concept
of embedding in $L_0$ was introduced and studied in \cite{KKYY}. As proved in \cite[Th. 6.3]{KKYY},
the unit balls of spaces that embed in $L_0$ are intersection bodies and also satisfy inequalities
(\ref{hyper-int}) and (\ref{arbmeas}).

\bigbreak
{\bf Acknowledgement.} I wish to thank the US National Science Foundation for support through 
grant DMS-1265155.

\end{document}